\input amstex
\expandafter\ifx\csname mntex.sty\endcsname\relax
  \expandafter\edef\csname mntex.sty\endcsname{%
    \catcode`\noexpand\@=\the\catcode`\@\relax}
\else \message{... already loaded}\endinput\fi

\catcode`\@=11

\def\W@#1{}\let\@message\message \def\message#1{}

\input amsppt.sty \relax

\pagewidth{29pc}

\let\message\@message \def\W@{\immediate\write\sixt@@n}

\def\styversion{2.2} \def\styname{MNTEX}
\W@{\styname.STY - Version \styversion}\W@{}

\outer\def\zh#1\endzh{\csname specialhead\endcsname#1\endspecialhead\nobreak}
\outer\def\ah#1\endah{\csname head\endcsname#1\endhead\nobreak}
\outer\def\bh#1\endbh{\csname subhead\endcsname#1\endsubhead\nobreak}
\outer\def\ch#1\endch{\csname subsubhead\endcsname#1\endsubsubhead\nobreak}

\outer\def\tm{\csname proclaim\endcsname}

\def\pf{\csname definition\endcsname}        

\def\demo{\csname definition\endcsname}
\def\enddemo{\enddefinition}

\def\de{\csname definition\endcsname}

\def\rk{\csname definition\endcsname}

\def\remark{\csname definition\endcsname}
\def\endremark{\enddefinition}

\def\ex{\csname definition\endcsname}

\def\example{\csname definition\endcsname}
\def\endexample{\enddefinition}

\def\translator{%
  \let\savedef@\translator 
  \def\translator##1\endtranslator{\let\translator\savedef@
    \edef\thetranslator@{\noexpand\nobreak\noexpand\medskip
      \noexpand\line{\noexpand\eightpoint\hfil
      \frills@{Translated by }{##1}}%
       \noexpand\nobreak}}%
  \nofrillscheck\translator}

\Monograph
\monograph@false
\let\headmark\eat@

\begingroup
\let\head\relax  
\let\ah\relax \let\bh\relax \let\ch\relax

\gdef\widestnumber{\begingroup
  \let\ah\relax \let\bh\relax \let\ch\relax
  \expandafter\endgroup\setwidest@}
\gdef\setwidest@#1#2{%
   \ifx#1\ah\setbox\tocheadbox@\hbox{#2.\enspace}%
   \else\ifx#1\bh\setbox\tocsubheadbox@\hbox{#2.\enspace}%
   \else\ifx#1\ch\setbox\tocsubheadbox@\hbox{#2.\enspace}%
   \else\ifx#1\key\refstyle A%
       \setboxz@h{\refsfont@\keyformat{#2}}%
       \refindentwd\wd\z@
   \else\ifx#1\no\refstyle C%
       \setboxz@h{\refsfont@\keyformat{#2}}%
       \refindentwd\wd\z@
   \else\ifx#1\page\setbox\z@\hbox{\quad\bf#2}%
       \pagenumwd\wd\z@
   \else\ifx#1\item
       \edef\next@{\the\revert@\rosteritemwd\the\rosteritemwd\relax
              \revert@{\the\revert@}}%
       \revert@\expandafter{\next@}%
       \setboxz@h{(#2)}\rosteritemwd\wdz@
   \else\message{\string\widestnumber\space not defined for this
      option (\string#1)}%
\fi\fi\fi\fi\fi\fi\fi}

\endgroup

\mathsurround=1.2pt
\def\zz{\leavevmode\kern-\mathsurround}

\def\punctcheck@{\futurelet\nextchar@\removemathsurround@}

\let\everymath@\everymath
\let\everydisplay@\everydisplay

\newtoks\everymath
\everymath@{\the\everymath}

\newtoks\everydisplay
\everydisplay@{\the\everydisplay}

\everymath{\aftergroup\punctcheck@}

\def\removemathsurround@{%
  {\edef\next{\noexpand\nextchar@}%
  \ifcat.\noexpand\nextchar@
    \expandafter\ifx\next,\else\expandafter\ifx\next.\else
      \expandafter\ifx\next;\else\kern-\mathsurround\fi\fi\fi
  \fi
  \let\nextchar@\relax
  }%
}

\let\@@eqno=\eqno

\let\@@leqno=\leqno

\def\eqno{\everymath{}\@@eqno}
\def\leqno{\everymath{}\@@leqno}

\def\mathcomma@{\mathpunct{\mkern1.2mu\mathchar"013B\mkern1.2mu}}%
\def\mathsemicolon@{\mathpunct{\mkern3mu\mathchar"003B\mkern3mu}}%
\def\mathcolon@{\mathrel{\mkern\thickmuskip\nonscript\mkern-\thickmuskip
  \mathchar"003A\mkern\thickmuskip\nonscript\mkern-\thickmuskip}}

\mathcode`\,="8000 \mathcode`\;="8000 \mathcode`\:="8000

\begingroup
\catcode`\,=\active \catcode`\;=\active \catcode`\:=\active
\global\everymath@{%
  \let,\mathcomma@ \let;\mathsemicolon@ \let:\mathcolon@
  \the\everymath}
\global\everydisplay@{%
  \let,\mathcomma@ \let;\mathsemicolon@ \let:\mathcolon@
  \def\.{\thinspace.}%
  \the\everydisplay}
\endgroup

\def\cprime{\begingroup\everymath{}\m@th$'$\endgroup}
\let\mz\cprime 

\font@\ninerm=cmr9
\font@\ninei=cmmi9    \skewchar\ninei='177
\font@\ninesy=cmsy9   \skewchar\ninesy='60
\font@\nineex=cmex9
\font@\ninebf=cmbx9
\font@\nineit=cmti9
\font@\ninesmc=cmcsc9
\font@\ninesl=cmsl9

\font@\ninemsa=msam9
\font@\ninemsb=msbm9
\font@\nineeufm=eufm9

\newtoks\ninepoint@
\def\ninepoint{\normalbaselineskip11\p@
 \abovedisplayskip10\p@ plus2.4\p@ minus4.2\p@
 \belowdisplayskip\abovedisplayskip
 \abovedisplayshortskip\z@ plus2.4\p@
 \belowdisplayshortskip5.6\p@ plus2.4\p@ minus3.2\p@
 \textonlyfont@\rm\ninerm \textonlyfont@\it\nineit
 \textonlyfont@\bf\ninebf \textonlyfont@\smc\ninesmc
 \textonlyfont@\sl\ninesl
 \ifsyntax@\def\big##1{{\hbox{$\left##1\right.$}}}%
  \let\Big\big \let\bigg\big \let\Bigg\big
 \else
  \textfont\z@\ninerm \scriptfont\z@\sevenrm
       \scriptscriptfont\z@\fiverm
  \textfont\@ne\ninei \scriptfont\@ne\seveni
       \scriptscriptfont\@ne\fivei
  \textfont\tw@\ninesy \scriptfont\tw@\sevensy
       \scriptscriptfont\tw@\fivesy
  \textfont\thr@@\nineex \scriptfont\thr@@\sevenex
   \scriptscriptfont\thr@@\sevenex
  \textfont\itfam\nineit \scriptfont\itfam\sevenit
   \scriptscriptfont\itfam\sevenit
  \textfont\bffam\ninebf \scriptfont\bffam\sevenbf
   \scriptscriptfont\bffam\fivebf
  \textfont\msafam\ninemsa \scriptfont\msafam\sevenmsa
   \scriptscriptfont\msafam\fivemsa
  \textfont\msbfam\ninemsb \scriptfont\msbfam\sevenmsb
   \scriptscriptfont\msbfam\fivemsb
  \textfont\eufmfam\nineeufm \scriptfont\eufmfam\seveneufm
   \scriptscriptfont\eufmfam\fiveeufm
 \setbox\strutbox\hbox{\vrule height7.5\p@ depth3.5\p@ width\z@}%
 \setbox\strutbox@\hbox{\raise.5\normallineskiplimit\vbox{%
   \kern-\normallineskiplimit\copy\strutbox}}%
 \setbox\z@\vbox{\hbox{$($}\kern\z@}\bigsize@1.2\ht\z@
 \fi
 \normalbaselines\ninerm\dotsspace@1.5mu\ex@.2326ex\jot3\ex@
 \the\ninepoint@}

\font@\twelverm=cmr10 scaled\magstep1
\font@\twelvei=cmmi10 scaled\magstep1    \skewchar\twelvei='177
\font@\twelvesy=cmsy10 scaled\magstep1   \skewchar\twelvesy='60
\font@\twelveex=cmex10 scaled\magstep1
\font@\twelvebf=cmbx10 scaled\magstep1
\font@\twelveit=cmti10 scaled\magstep1

\font@\twelvemsa=msam10 scaled\magstep1
\font@\twelvemsb=msbm10 scaled\magstep1
\font@\twelveeufm=eufm10 scaled\magstep1

\newtoks\twelvepoint@
\def\twelvepoint{\normalbaselineskip15\p@
 \abovedisplayskip10\p@ plus2.4\p@ minus7.2\p@
 \belowdisplayskip\abovedisplayskip
 \abovedisplayshortskip\z@ plus2.4\p@
 \belowdisplayshortskip5.6\p@ plus2.4\p@ minus3.2\p@
 \textonlyfont@\rm\twelverm \textonlyfont@\bf\twelvebf
 \textonlyfont@\it\twelveit
 \ifsyntax@\def\big##1{{\hbox{$\left##1\right.$}}}%
  \let\Big\big \let\bigg\big \let\Bigg\big
 \else
  \textfont\z@\twelverm \scriptfont\z@\tenrm
       \scriptscriptfont\z@\sevenrm
  \textfont\@ne\twelvei \scriptfont\@ne\teni
       \scriptscriptfont\@ne\seveni
  \textfont\tw@\twelvesy \scriptfont\tw@\tensy
       \scriptscriptfont\tw@\sevensy
  \textfont\thr@@\twelveex \scriptfont\thr@@\tenex
   \scriptscriptfont\thr@@\sevenex
  \textfont\itfam\twelveit \scriptfont\itfam\tenit
   \scriptscriptfont\itfam\sevenit
  \textfont\bffam\twelvebf \scriptfont\bffam\tenbf
   \scriptscriptfont\bffam\sevenbf
  \textfont\msafam\twelvemsa \scriptfont\msafam\tenmsa
   \scriptscriptfont\msafam\sevenmsa
  \textfont\msbfam\twelvemsb \scriptfont\msbfam\tenmsb
   \scriptscriptfont\msbfam\sevenmsb
  \textfont\eufmfam\twelveeufm \scriptfont\eufmfam\teneufm
   \scriptscriptfont\eufmfam\seveneufm
 \setbox\strutbox\hbox{\vrule height12\p@ depth5\p@ width\z@}%
 \setbox\strutbox@\hbox{\raise.5\normallineskiplimit\vbox{%
   \kern-\normallineskiplimit\copy\strutbox}}%
 \setbox\z@\vbox{\hbox{$($}\kern\z@}\bigsize@1.2\ht\z@
 \fi
 \normalbaselines\twelverm\dotsspace@1.5mu\ex@.2326ex\jot3\ex@
 \the\twelvepoint@}

\newbox\plenumcopyright

\output={\faaoutput}
\def\faaoutput{\shipout\vbox{\pagebody\makefootline}%
  \advancepageno
  \ifnum\outputpenalty>-\@MM \else\dosupereject\fi}
\def\pagebody{\advance\vsize-0.6\ht\plenumcopyright
\vbox to\vsize{\boxmaxdepth\maxdepth \pagecontents}}
\def\makeheadline{\vbox to\z@{\vskip-22.5\p@
  \line{\vbox to8.5\p@{}\the\headline}\vss}\nointerlineskip}
\def\makefootline{\baselineskip24\p@\line{\the\footline}}
\def\dosupereject{\ifnum\insertpenalties>\z@ 
  \line{}\kern-\topskip\nobreak\vfill\supereject\fi}

\def\pagecontents{\ifvoid\topins\else\unvbox\topins\fi
  \dimen@=\dp\@cclv \unvbox\@cclv 
  \ifvoid\footins\else 
    \vskip\skip\footins
    \footnoterule
    \unvbox\footins\fi
      \ifr@ggedbottom \kern-\dimen@ \vfil \fi}
\def\footnoterule{\vskip2pt plus 1pt
  \hrule width 2in \kern 2.6\p@} 
\footline{\ifvoid\plenumcopyright\sevenrm\ifodd\pageno
  \hss\the\pageno\else\the\pageno\hss\fi\else\box\plenumcopyright\fi}
\outer\def\enddocument{\par
  \add@missing\endRefs
  \add@missing\endroster \add@missing\endproclaim
  \add@missing\enddefinition
  \add@missing\enddemo \add@missing\endremark \add@missing\endexample
 \ifmonograph@ 
 \else
 \nobreak
 \thetranslator@
 \count@\z@ \loop\ifnum\count@<\addresscount@\advance\count@\@ne
 \csname address\number\count@\endcsname
 \csname email\number\count@\endcsname
 \repeat
\fi
 \vfill\supereject\end}

\csname mntex.sty\endcsname
\pagewidth{38.33pc} \pageheight{52.5pc}
\hoffset=-1.25cm
\voffset=-27pt
\baselineskip=12pt plus 0.2pt minus 0.2pt
\abovedisplayskip=7pt plus 2pt minus 2pt
\belowdisplayskip=7pt plus 2pt minus 2pt
\abovedisplayshortskip=5pt plus 1pt minus 3pt
\belowdisplayshortskip=5pt plus 1pt minus 3pt

\def\titleinfo#1{\flushpar\eightit Mathematical Notes,
Vol.~\Russvol, No.~\Russno,~19\Russyr
   \vglue 2.5pc \flushpar{\twelvepoint \bf #1}\rm}

\def\extratitleline#1{\smallskip\flushpar{\twelvepoint\bf #1}\rm}

\def\authorinfo#1#2#3{\medskip\medskip \bf#1
\hfill\rm #2\medskip\medskip
\setbox\plenumcopyright=\vbox{\vskip12pt\hrule width \hsize \kern 3pt
\baselineskip9.5pt\eightrm Translated from {\eightit Matematicheskie
Zametki}, Vol.~\Russvol, No.~\Russno, pp.~\Russpages, \Russmonth, 19\Russyr.
\newline\indent #3.\newline
\centerline{\vrule height12pt width0pt\ifodd\pageno{}\else{\sevenrm\the\pageno}\fi
\hfill 0001--4346/\Russyr/\Russvol\Englno--{\ifnum \pageno<10 000\the\pageno%
         \else\ifnum\pageno<100 00\the\pageno
         \else\ifnum\pageno<1000 0\the\pageno
         \else\ifnum\pageno<9999 \the\pageno\fi\fi\fi\fi}\,\$\copyfee\quad
\copyright19\Englyr\ Plenum
Publishing Corporation\hfill\ifodd\pageno{\sevenrm\the\pageno}\else{}\fi}}
}

\def\authorinfopage#1#2#3{\medskip\medskip \bf#1
\hfill\rm #2\medskip\medskip
\setbox\plenumcopyright=\vbox{\vskip12pt\hrule width \hsize \kern 3pt
\baselineskip9.5pt\eightrm Translated from {\eightit Matematicheskie
Zametki}, Vol.~\Russvol, No.~\Russno, p.~\Russpage, \Russmonth, 19\Russyr.
\newline\indent #3.\newline
\centerline{\vrule height12pt width0pt\ifodd\pageno{}\else{\sevenrm\the\pageno}\fi
\hfill 0001--4346/\Russyr/\Russvol\Englno--{\ifnum \pageno<10 000\the\pageno%
         \else\ifnum\pageno<100 00\the\pageno
         \else\ifnum\pageno<1000 0\the\pageno
         \else\ifnum\pageno<9999 \the\pageno\fi\fi\fi\fi}\,\$\copyfee\quad
\copyright19\Englyr\ Plenum
Publishing Corporation\hfill\ifodd\pageno{\sevenrm\the\pageno}\else{}\fi}}
}

\redefine\Refs{\bgroup\vskip10pt%
\centerline{\bf References}\nobreak\vskip10pt\leftskip5pt\frenchspacing
\eightpoint\nobreak}
\redefine\endRefs{\egroup}
\def\authoraddress#1{\vskip10pt minus 1pt
\leftline{\eightpoint\smc\ \ \ \ #1}}

\def\transl#1{\vskip10pt minus 1pt \rightline{\eightpoint Translated by #1}}
\def\ed--transl#1{\vskip10pt minus 1pt
\eightbf EDITOR: \
\underbar{$\hphantom{\text{Alex-Alexandrovskii}}$}\hfill
\eightrm Translated by #1}
\def\bcskip{\vglue30pt\hrule width \hsize
\kern3pt\vglue30pt}
\TagsOnRight
\NoRunningHeads
\def\bS{\pmb\S}
\define\q{\quad}
\define\qq{\qquad}

\define\nlb{\nolinebreak}

\define\hph#1{\hphantom{#1}}
\comment

\def\sba#1{_{\lower1pt\hbox{$\scriptstyle#1$}}}     
\def\sbb#1{_{\lower2pt\hbox{$\scriptstyle#1$}}}
\def\spa#1{^{\raise1pt\hbox{$\scriptstyle#1$}}}     
\def\spb#1{^{\raise2pt\hbox{$\scriptstyle#1$}}}
\define\z#1{{\bold#1}}          
\define\zk#1{\boldkey#1}        
\define\zs#1{\boldsymbol#1}     
\endcomment


\def\script{\Cal}
\def\cal{\Cal}
.pk scaled 1000
.pk scaled 1000
.pk scaled 1000
.pk scaled\magstep1



\pageno=510
\scrollmode

\define\endbad{\egroup}
\define\mes{\operatorname{mes}}
\define\Int{\operatorname{int}}
\define\R{\Bbb R}
\define\FF{\script F}
\define\dom{\operatorname{dom}}
\let\hlop\!
\let\point.

\let\phi\varphi
\let\epsilon\varepsilon
\let\cal\script

\def\Russvol{61}%
\def\Russno{4}%
\def\Englno{34}%
\def\Russpages{612--622}
\def\Russmonth{April}%
\def\Russyr{97}%
\def\Englyr{97}%
\def\copyfee{18.00}%

\titleinfo
{The Quantum Stochastic Differential Equation Is Unitarily Equivalent}
\extratitleline{to a Symmetric
Boundary Value Problem for the Schr\"odinger Equation}
\authorinfo
{A.~M.~Chebotarev}
{UDC 519.217}
{Original article submitted December 11, 1996}

\topmatter
\abstract\nofrills
{\smc Abstract}.
We prove that the solution of the Hudson--Parthasarathy
quantum stochastic  differential equation in the Fock space coincides with
the solution of a symmetric
boundary value problem for the Schr\"odinger equation in the interaction
representation generated by the energy  operator of the environment.
The boundary conditions describe the jumps in the phase and the amplitude
of the Fourier transforms of the Fock vector  components
as  any of its arguments changes the sign.
The corresponding Markov evolution  equation (the Lindblad equation
 or the ``master equation'') is derived from the
boundary value problem for the Schr\"odinger equation.
\newline\newline\noindent
{\smc Key words}: \,     
singular perturbation, quantum noises, boundary value problem in the Fock
space, open quantum system, master equation, quantum dynamical semigroup.
\endabstract
\endtopmatter

\document

\head\bS1. Introduction \endhead

Our interest in the theory of
quantum stochastic differential equations is
primarily due to the fact that it can be used  
to represent  solutions of the
Markov evolution equation,  which generalizes the Heisenberg, 
 Kolmogorov--Feller, and  heat equations \cite{1--4}.

Solutions of Markov evolution  equations are studied in the theory of
open quantum systems \cite{5--8}; they are called
{\it quantum dynamical semigroups}.
One of the consequences of the stochastic representation of quantum
dynamical semigroups in  von Neumann algebras is a simple relationship
between the conservativity  of a quantum dynamical semigroup and the
isometric property of a solution of the quantum stochastic differential 
equation.
The conservativity (or the unital property) of a quantum dynamical
semigroup means
that the unit of the  operator algebra is preserved; this corresponds
to the preservation of the total probability by the solutions of
the
equations for the transition probabilities of the stochastic processes in
the classical and quantum cases.
Violations of this property are related to critical phenomena such as
explosions, escape of the solution to infinity,
accumulation of infinitely many  discontinuities, 
nonunique  solvability of the Cauchy problem,
the existence of growing solutions to formally dissipative
equations, etc. \cite{9--10}.
The study of conditions necessary and sufficient for 
the existence or nonexistence
of such phenomena is currently of considerable interest~\cite{11--18}.

Sufficient criteria for   solutions of
quantum stochastic differential equations 
to be isometric or unitary
were obtained  in recent years 
mainly by
perturbation methods~\cite{19--21}.
However, symmetric operators  responsible for
these
properties of  solutions were not found.
The difficulties encountered in this approach 
include the violation of the
group property by solutions of 
quantum stochastic differential equations  \cite{22} and the
fact that their formal generators fail to be symmetric
on domains consisting of smooth functions.
More precisely \cite{23}, the formal generators for the 
Schr\"odinger equations  unitarily equivalent to
quantum stochastic differential equations  have the form of
dissipative operators perturbed by singular
quadratic forms  in the sense of \cite{24--25}.
Since the dissipative part of the 
generator has an  unbounded anti-Hermitian part, the general
theory of 
self-adjoint extensions of  symmetric
operators bounded below and perturbed by singular quadratic
forms \cite{24--26} does not apply directly to this case.

In the present paper, we consider quantum stochastic
differential equations whose solutions  in the interaction
representation  are limits of solutions of a family of
Schr\"odinger  equations in the Fock space \cite{23}.
The fundamental property of functions belonging to the range of
the resolvent of the limit unitary group is 
that they experience  {\it amplitude} and {\it
phase jumps} at the points where
the quadratic form of \linebreak

\pagebreak\noindent
the generator has singularities.
Evaluating the action of the generator of the limit group on 
functions with amplitude and phase  jumps,  we obtain a symmetric
operator; its self-adjointness can be verified independently
under assumptions less restrictive than those used to construct
the resolvent explicitly.
Thus, we prove the equivalence of the quantum stochastic
differential equation and 
the boundary value 
problem for  the Schr\"odinger equation in the Fock space.
Straightforward computations permit one to derive
the Markov evolution  equation  directly from the boundary value
problem for 
the Schr\"odinger equation.

We start from a simple example showing  the
distinction between the
formal weak limit, the strong resolvent limit of a sequence of
generators of 
unitary groups, and the action of the generator of the limit unitary
group on the range of its resolvent.

\head\bS2. The weak limit and the resolvent limit of a family
of generators of one-parameter unitary groups 
\endhead

Consider the  one-parameter family of unitary
groups $\exp\{it\widehat H_\alpha\} = U_t^{(\alpha)} $ in
$L_2(\R)$ given by
$$
U_t^{(\alpha)}\psi(x)
=\psi (x-t) \exp\biggl\{i\lambda\int_0^td\tau\, V_\alpha (x-t + \tau)
\biggr\}, \qquad
x,\lambda\in\R,
$$
where $V_\alpha (x) = (2\pi\alpha) ^{-1 / 2} \exp\{-x^2 / 2\alpha\} $,
$\alpha\in\R_+$.
Obviously, $V_\alpha (x) \to\delta (x) $ as $\alpha\to +
0$.
Therefore, the {\it weak limit\/} of the family of essentially self-adjoint
operators
$\widehat H_\alpha=i\partial_x + \lambda V_\alpha (x) $ is described by the
bilinear form
$$
\widehat H_w [ \varphi,\psi ] = (\varphi,\widehat H_w\psi) =i (
\varphi,\psi') + \lambda\overline\varphi (0) \psi (0),
\tag2.1
$$
which is well defined on  $W_2^1 (\R) $.
On the other hand,
$$
\int_0^td\tau\, V_\alpha (x-t + \tau) \to I_{[ 0, t)} (x),
$$
where $I_T (x) $ is the characteristic function of a Borel set
$T\subseteq\R$.
Hence, the {\it strong limit\/} of the family $U_t^{(\alpha)}$
is 
$$
\lim_{\alpha\to+0}U_t^{(\alpha)}\psi(x)
=U_t\psi(x)
=e^{it\widehat {\bold H}_R} \psi (x) =\psi (x-t) e^{i\lambda I_{[
0, t)} (x)}.
$$
Note the identity $e^{i\lambda I_{[ 0, t)} (x)} = (e^{i\lambda} -1)
I_{[ 0, t)} (x) + 1$.
Therefore, the bilinear form of the {\it resolvent limit\/}
$\widehat {\bold H}
_R=\text{$r$-lim}\,\widehat H_\alpha$ is well defined on $W_2^1 (\R) $ 
and is given by
$$
\widehat {\bold H}_R [ \varphi,\psi ] =i^{-1} \lim_{t\to0} \frac d {
dt} (\varphi, U_t\psi) \Big|_{t=0} =i (\varphi,\psi') + i (e^{
i\lambda} -1) \overline\varphi (0) \psi (0).
\tag2.2
$$
Comparing~\thetag{2.1} with~\thetag{2.2}, we conclude that $\widehat {\bold
H}_R=\text{$r$-lim}\,\widehat H_\alpha \ne \text{$w$-lim}\, H_\alpha=
\widehat H_w$.

The range of the resolvent of the limit unitary group $U_t$ is the natural
domain of the generator $\widehat {\bold H}_R$ and 
can be described explicitly in our example:
$$
R_\mu\psi(x)
=\int_0^\infty dt\, e^{-\mu t} \psi (x-t) + \theta (x) (e^{i\lambda}
-1) e^{-\mu x} \int_0^\infty dt\, e^{-\mu t} \psi (-t),
$$
where $\theta (x) $~ is the characteristic function of the semiaxis $\R_+ $.
This structure of the resolvent means that the functions
belonging to the domain of
$\widehat {\bold H}_R$ have
a phase jump at the origin $x=0$:
$$
\lim_{x\to+0}R_\mu\psi(x)
=e^{i\lambda}\lim_{x\to-0}R_\mu\psi(x).
$$
Therefore, the domain $\cal D_\lambda$ of the generator
$\widehat {\bold H}_R$ consists of the functions
$$
\psi: \quad \psi\in W_2^1 (\R\setminus\{0\}), \qquad
\lim_{x\to+0}\psi(x)
=e^{i\lambda}\lim_{x\to-0}\psi(x).
\tag2.3
$$
The operator $\widehat {\bold H}_R$ acts as
$i\partial_x$ for $x\ne0$.
For functions $\psi\in \cal D_\lambda$, the left and the right
limits    at the 
origin exist
by virtue of the embedding $W_2^1 (\R\setminus\{0\})
\subset C (\R\setminus\{0\}) $.
The symmetry property of  $\widehat {\bold H}_R$ follows
>from  integration by parts and from the identity
$$
\overline {\phi (x)} \psi (x) \big|_{+ 0} ^{-0} =0,
$$
which holds for $\varphi,\psi\in \cal D_\lambda$;
the self-adjointness follows  \cite{27} from the fact that 
for $\mu=\pm1$
the
equation $ (\widehat {\bold H}_R + i\mu) \psi
(x) =f (x) $, $x\ne0$, with the boundary condition~\thetag{2.3
} is solvable in $\cal D_\lambda$ for every right-hand side
$f\in L_2 (\R) $.

\head\bS3. The weak and the strong resolvent limit of solutions
of the Schr\"odinger equation in the Fock space \endhead

Suppose that  $\cal H$ is a Hilbert space,
$\Gamma^S (L_2 (\R)) $ is the symmetric Fock space and $\bold h=\cal
H\otimes\Gamma^S (L_2 (\R)) $ is their tensor product.
For $v, g\in L_2 (\R) $,
by $A^+ (v)$ and $ A (g) $ we denote the standard creation and annihilation
operators in $\Gamma^S (L_2 (\R)) $. Consider the family of
Schr\"odinger equations $\partial_t\psi_t=iH\psi_t$
with  self-adjoint   Hamiltonian $\bold H=H_0\otimes I +
I\otimes\widehat E + H_{\Int} $ in $\bold h$:
$$
H_{\Int}
=K\otimes A^+ (g) A (g)
+ R\otimes A^+ (f) + R^* \otimes A (f), \qquad \widehat E=\int\omega a^+ (
\omega) a (\omega) \, d\omega.
\tag3.1
$$
>From now on, it is assumed for simplicity   that the operators $H_0,K,$ and
$R\in\cal C (\cal H) $ commute and have a joint spectral family $E_\lambda$.
More precisely, let
$$
H_0=\int\nu (\lambda) \, dE_\lambda, \qquad K=\int\lambda\, dE_\lambda,
\qquad R=\int\rho (\lambda) e^{i\Phi (\lambda)} \, dE_\lambda,
$$
where $\nu,\rho,$ and $\Phi$ are measurable real functions corresponding
to the operators $H_0$, $K$, and $R$, so that  the operators
$H_0$ and $K$ are  self-adjoint and $R$ is normal.

We denote by $\widehat P_t (\lambda) $ the one-parameter group of
unitary operators in $L_2 (\R) $ with the generator $\widehat N (\lambda)
=\omega + \lambda| g\rangle\langle g| $.
The unitary group $U_t=\exp\{i\bold Ht\} $ can be constructed
by  transforming  the Hamiltonian to the canonical form
\cite{28}, or by  using the interaction representation
generated by the operator $I\otimes\widehat E +
K\otimes A^+ (g) A (g) $ \cite{23}.
Theorem 3.1 describes the action of $U_t$ on coherent vectors.

\proclaim{Theorem 3.1}
The unitary one-parameter group  $U_t=\exp\{it\bold H\} $, where
$\bold H=H_0\otimes I + I\otimes\widehat {\bold E} + H_{\Int} $
and  $H_{\Int} $ is defined by Eq.~\thetag{3.1},
acts as follows\rom:
$$
U_t h\otimes\psi (v) =\int e^{i\nu (\lambda) t} \, dE_\lambda
h\otimes\psi\bigl(v_t(\lambda)\bigr)
\exp\biggl\{i\rho (\lambda) e^{-i\Phi (\lambda)} \int_0^t\bigl (f,
v_s (\lambda) \bigr) \, ds\biggr\},
\tag3.2
$$
where
$$
v_t(\lambda)
=\widehat P_t (\lambda) v + i\rho (\lambda) e^{i\Phi (\lambda)}
\int_0^t\widehat P_s (\lambda) f\, ds.
$$
\endproclaim

Let $L_{2,\widetilde1} ^+ (\R)\subset \widetilde W_2^1 (\R) $ 
be the set of functions with 
positive  absolutely integrable Fourier transform; 
$\widetilde W_2^1 (\R) $ is the
Fourier transform of the Sobolev space $W_2^1 (\R) $.
Let $f, g\in L_{2,\widetilde1} ^+ (\R) $ be real functions
such that $f (0) =g (0) = 1 / \sqrt {2\pi} $, and let 
$f^{(\alpha)} (\omega)=f (\alpha\omega) $,
$g^{(\alpha)}(\omega)=g(\alpha\omega)$.

\proclaim{Lemma 3.1 {\rm (about four limits)}} Let $\widehat P^{(
\alpha)}_t (\lambda) $ be the one-parameter unitary group
in $L_2 (\R) $ with  generator $\widehat N_\alpha (\lambda) =\omega +
\lambda| g^{(\alpha)} \rangle\langle g^{(\alpha)}| $, and let
$\widehat\pi_T=\FF_{t\to\omega} I_T (t) \FF^*_{\omega\to t} $ be
a
family of projections in $\cal H (L_2 (\R)) $.
Then  the following  limits exist  as $\alpha\to0$\rom:
\roster
\item"\rom{1)}"
$\dsize\lim\int^t_0ds\,
\bigl (g^{(\alpha)},\widehat P^{(\alpha)}_s (\lambda) f^{(\alpha)} \bigr)
=(2-i\lambda)^{-1}$\rom;
\item"\rom{2)}"
$\dsize \text{\rm$w$-lim}\,\int^t_0 ds\,
\widehat P^{(\alpha)}_s (\lambda) f^{(\alpha)} (\omega) =e^{
i\omega t} \widetilde I_{[0, t)} (\omega) (1-i\lambda / 2) ^{-1}
$\rom;
\item"\rom{3)}"
$\text{\rm$w$-lim}\,\bigl (g^{(\alpha)},\widehat P^{(\alpha)}_t (\lambda)
v\bigr) = (1-i\lambda / 2) ^{-1} \FF^*_{\omega\to t} v$\rom;
\item"\rom{4)}"
$\text{\rm$s$-lim}\,\widehat P^{(\alpha)}_t (\lambda)
=\exp\{iZ(\lambda)\widehat\pi_{(0,t)}\}
=\widehat P_t (\lambda) $,
$\exp\{iZ (\lambda) \} = (2 + i\lambda) / (2-i\lambda) $.
\endroster
\endproclaim

Let us consider the family $\bold H_\alpha$ of self-adjoint Hamiltonians
\thetag{3.1} parametrized by the functional arguments
$f^{(\alpha)} (\omega) $, $g^{(\alpha)} (\omega) $ 
of the creation and annihilation operators.
Lemma~3.1 in~\cite{23} justifies the passage to the limit  in
Eq.~\thetag{3.2} as $\alpha\to + 0$.
The limits~1)--4) correspond to the passage to the limit  in the
four components 
of the solution~\thetag{3.2} depending on $\alpha$.
By substituting~1)--4) into Eq.~\thetag{3.2}, we obtain the
limit unitary group 
$U_t=\exp\{i\bold Ht\} =\text{$s$-lim}_{\alpha\to0} \exp\{i\bold
H_\alpha t\}$: 
$$
\align
U_th\otimes\psi(v)
&=\int e^{iH(\lambda)t}\,dE_\lambda h
\psi\Bigl(e^{iZ(\lambda)\widehat\pi_{[0,t)}}e^{i\omega t}v
+\frac{2i}{2-i\lambda}\rho(\lambda)e^{i\Phi(\lambda)}
\widetilde I_{[0,t)}\Bigr)
\\ &\qquad\times
\exp\Bigl\{\rho(\lambda)e^{-i\Phi(\lambda)}
\frac{2i}{2-i\lambda}(\widetilde I_{[0,t)},e^{i\omega t}v)\Bigr\},
\tag3.3
\endalign
$$
where $H (\lambda) =\nu (\lambda) + i\rho^2 / (2-i\lambda) $.
Let  $W$ and $L$ be
the operators with spectral densities $(2+i\lambda)/(2-i\lambda)$
and $2 i / (2-i\lambda) \,\rho (\lambda) e^{i\Phi (\lambda)} $,
respectively. 
Then $ {2 i} \rho (\lambda) e^{-i\Phi (\lambda)} / (2-i\lambda) $
is the spectral density of the operator $-L^* W$, where
$W= (2 + iK)/ (2-iK) $
is the Cayley transformation  of the self-adjoint operator $2
K$, and 
$L=2 i / (2-iK) \, R$ is a densely defined operator such
that $\dom L\supseteq\dom
R$, and $H (\lambda) $ is the spectral function of the operator
$$
iG=H_0-\frac14 L^* KL + \frac i2 L^* L.
$$
In the following, we suppose that the operator $-G$ is the generator
of a strongly continuous one-parameter contraction semigroup
$W_t=\exp\{-Gt\}$ 
in $\cal H$, and moreover,
$$
D=\dom H\cap\dom L^* L\subseteq\dom G\subseteq\dom L, \qquad G^* \phi +
G\phi=L^* L\phi \quad \forall\phi\in D,
$$
where the operator $H=-H_0 + \frac14 L^* KL$ is symmetric on $D$.
In this notation, the bilinear form
$$
\bold H [ g\otimes\psi (f), h\otimes\psi (v) ]
=\lim_{t\to+0}
\frac1 i\,\frac d {dt} \bigl (g\otimes\psi (f), U_t\, h\otimes\psi (v)
\bigr),
$$
where $g, h\in D\subseteq\cal H$ and $f, v\in\widetilde W_2^1 (\R) $, acts
as follows:
$$
\align
\bold H [ g\otimes\psi (f), h\otimes\psi (v) ]
&=e^ {(f, v)_{L_2}} \biggl (i (g, Gh)_{\cal H} -i (g, Lh)_{
\cal H} \overline {\widetilde f (0)} + i (Lg, Wh)_{\cal H} \widetilde
v (0) \\
&\quad + (g, h)_{\cal H} \int d\omega\, \overline f (\omega) g
(\omega) \omega + i (g, (I-W) h)_{\cal H} \overline {\widetilde f
(0)} \widetilde v (0) \biggr).
\tag3.4
\endalign
$$
The bilinear form~\thetag{3.4} has the regular dissipative
component $\frac i2 L^ * L\otimes I$ and 
the {\it singular\/} component~\cite{24--\nlb 25}
$$
-iL\otimes A^+ \biggl (\frac1 {\sqrt {2\pi}} \biggr) + iL^* W\otimes
A\biggl (\frac1 {\sqrt {2\pi}} \biggr) + i (I-W) \otimes A^+ \biggl (
\frac1 {\sqrt {2\pi}} \biggr)
A\biggl(\frac1{\sqrt{2\pi}}\biggr),
$$
vanishing on the subset $\cal D_0$
$$
\cal D_0=\{\Phi:\Phi=h\otimes\psi (v),
V\in\widetilde W_2^1 (\R), \widetilde v (0) =0, h\in D\},
$$
which is total in $\bold h$.
At the same time, the weak limit of the sequence $H_\alpha$ is
described  by a bilinear form with different
regular and singular components:
$$
\align
\lim_{\alpha\to0} {H_\alpha} [ g\otimes\psi (f), h\otimes\psi (v) ]
&=e^{(f, v)_{L_2}} \biggl((g, H_0 h)_{\cal H} + (g, h)_{\cal H}
\int\omega d\omega\,\overline f(\omega) g(\omega) \\
&\qquad + (g, Kh)_{\cal H} \overline {\widetilde f (0)}
\widetilde v (0) + (g, Rh)_{\cal H}
\overline {\widetilde f(0)} - (Rg, h)_{\cal H} {\widetilde v(0)} \biggr).
\endalign
$$
Just as in the example in \S2, here we have 
$$
r-\lim_{\alpha\to0} {\bold H}_\alpha \ne w-\lim_{\alpha\to0} {\bold H}_\alpha.
$$

Let us consider the set function $u(s,t)=J_sU_{t-s}J_t^*$,
where $J_t$ is the unitary one-parameter group generated by
the operator $\widehat E=\int
d\omega\,\omega a^+ (\omega) a (\omega) $ \cite{7}, which acts in
$\Gamma^S$ as $J_t\psi (v) =\psi (e^{i\omega t} v) $.
Taking into account the fact that $J^*_t\psi (v) =\psi (e^{-i\omega t} v) $, 
>from Eq.~\thetag{3.3} we obtain
$$
\align
u(T)h\otimes\psi(v)
=\int e^{iH (\lambda) \mes T} \, dE_\lambda\, h\otimes\psi\Bigl (e^{
iZ (\lambda) \widehat\pi_T} v
+i\rho(\lambda)e^{i\Phi(\lambda)}
\frac2 {2-i\lambda} \widetilde I_T\Bigr) \\
\qquad\times \exp\Bigl\{
i\rho (\lambda) e^{-i\Phi (\lambda)} \frac2 {2-i\lambda} (\widetilde
I_T, v) \Bigr\}.
\tag3.5
\endalign
$$
For disjoint sets, the family of commuting operators $u (T) $ satisfies
the composition law $u (T_1\cup T_2) =u (T_1) u (T_2) $, and the
differential of the bilinear form $ (g\otimes\psi (f), u (0, t) h\otimes\psi
(v)) $ satisfies the weak quantum stochastic differential equation \cite{2}
$$
d\bigl (h\otimes\psi (v), u (0, t) h\otimes\psi (v) \bigr) =i\bigl (
h\otimes\psi (v), u (0, t) H (dt_+) h\otimes\psi (v) \bigr),
\tag3.6
$$
where
$$
\gather
iH (T) =M (T) =\int_Tdt\, J_t (\bold H-\widehat E\otimes I) J^*_t,
\\
M (T) =-G\otimes\mes T + L\otimes A^+ (T) -L^* W\otimes A (T) + (W-I)
\otimes\Lambda (T).
\endgather
$$
Thus, the following statement holds.

\proclaim{Theorem 3.2}
The family of  solutions of the
Schr\"odinger equation 
with  Hamiltonian~\thetag{3.1} in $\bold h$
strongly converges,
up to the unitary transformation $J^*_t$,
to the solution $u (0,t) $ of the stochastic
equation~\thetag{3.6}, $u (0, t) =\text{$s$-lim}_{\alpha\to0}
U^{(\alpha)}_tJ^*_t$.  
\endproclaim

In what follows, 
we shall see that, quite unexpectedly, the natural domain of the
generator of the group~$U_t$ 
(the range of the resolvent) does not contain functions
on which the bilinear form~\thetag{3.4} is well defined:
the Fourier transforms of the  Fock components of the resolvent
experience 
amplitude and phase jumps at the points where their arguments
change sign.

\head\bS4. Surprising properties of the resolvent \endhead

Consider  the following vector
$\Phi\in\bold h=\cal H\otimes\Gamma^S\bigl (L_2 (\R) \bigr) $:
$$
\Phi=R_\mu h\otimes\psi (v) =\int_0^\infty dt\, e^{-\mu t} U_th\otimes\psi
(v)
=\{\Phi_n(\omega)\},
\qq
\Phi_n (\,\cdot\,) \: \R^n\to\cal H, \q
\omega=\{\omega_1,\dots,\omega_n\},
$$
with components~\thetag{3.3}
$$
\gather
\Phi_n(\omega)
=\int_0^\infty dt\,\exp\biggl\{-\mu (G + t -L^* W \int_0^t\widetilde v (
-\tau) \, d\tau\biggr\}
\phi_{n,t}(\omega),
\\
\phi_{n,t}(\omega)
=\prod_1^n\bigl ((W-1) \pi_{[ 0, t)} e^{i\omega_k t} v (\omega_k) +
e^{i\omega_k t} v (\omega_k) + L\widetilde I_{[ 0, t)} (\omega_k)
\bigr) h,
\endgather
$$
where $L$, $W$, and $G$ are the above-described commuting operators 
with spectral densities
$$
L(\lambda)
=2 i\rho (\lambda) e^{-i\Phi (\lambda)} (2-i\lambda) ^{-1}, \qquad
W(\lambda)
=e^{iZ (\lambda)}, \qquad
G(\lambda)
=-i\nu (\lambda) + \frac {\rho (\lambda) ^2} {(2-i\lambda)}.
$$
Let  $\widetilde\phi_{n, t}(\tau) $ be the Fourier transform
of  $\phi_{n, t}(\omega)$,
$$
\widetilde\phi_{n,t}(\tau)
=\prod_1^n\bigl ((W-I) I_{[ 0, t)} (\tau_k) \widetilde v (\tau_k-t) +
\widetilde v (\tau_k-t) + LI_{[ 0, t)} (\tau_k) \bigr) h,
$$
where $\tau=\{\tau_1,\dots,\tau_n\} $.
Let $\cal K\subset\{1,\dots, n\} $, and let $\cal K^a$  be the
complement of $\Cal K$. Set
$$
P_{\cal K, t} ^{(n)} (\tau)
=\prod_{k\in\cal K} \bigl ((W-I) \widetilde v (\tau_k-t) + L\bigr) I_
{[ 0, t)} (\tau_k) \in\cal B (\bold h).
$$
Then
$$
\widetilde\phi_{n,t}(\tau)
=\sum_{\cal K} \biggl (P_{\cal K, t} ^{(n)} (\tau) \prod_{
m\in\cal K^a} \widetilde v (\tau_m-t) \biggr) h.
\tag4.1
$$

Obviously, the functions $P_{\cal K, t} ^{(n)} (\tau) $ have
discontinuities of the first kind at the points where the variables $\tau_k$
change the sign:
$$
\aligned
\lim_{\tau_k\to-0} P_{\cal K, t} ^{(n)} (\tau)& =I_{\cal K^a}
(k) P_{\cal K, t} ^{(n)} (\tau),
\\
\lim_{\tau_k\to + 0} P_{\cal K, t} ^{(n)} (\tau) &=I_{\cal K^a}
(k) P_{\cal K, t} ^{(n)} (\tau) + \bigl ((W-I) \widetilde v (-t) +
L\bigr) P_{\cal K\setminus\{k\}, t} ^{(n-1)} (\tau).
\endaligned
\tag4.2
$$
Therefore, Eq.~\thetag{4.2} implies
$$
P_{\cal K, t} ^{(n)} (\tau) \big|_{\tau_k= + 0} ^{\tau_k=-0}
=-\bigl ((W-I) \widetilde v (-t) + L\bigr) P_{\cal K\setminus\{k\},
t} ^{(n-1)} (\tau) I_{\cal K} (k).
\tag4.3
$$
Let us calculate the jump of $\widetilde\phi_{n, t} (\tau) $
at the points where $\tau_k$ changes sign.
Note that
$$
\lim_{\tau_k\to-0} \widetilde\phi_{n, t} (\tau) =\widetilde v (-t)
\widetilde\phi_{n-1, t} (\tau_1,\dots,\tau_{k-1},\tau_{k + 1}
,\dots,\tau_n).
$$
Taking into account Eqs.~\thetag{4.1} and~\thetag{4.3}, we find
the amplitude and phase jumps  of the functions from the
domain of the infinitesimal operator of~$U_t:$
$$
\lim_{\tau_k\to + 0} \widetilde\phi_{n, t} (\tau) =W\lim_{
\tau_k\to-0} \widetilde\phi_{n, t} (\tau)
+L\widetilde\phi_{n-1,t}
(\tau_1,\dots,\tau_{k-1},\tau_{k + 1},\dots,\tau_n).
\tag4.4
$$
Let $\cal D_{W, L} =D\otimes\Gamma^S (\widetilde W_2^1 (
\R\setminus\{0\}))\subset \bold h $ be the vector subspace of
elements satisfying  condition~\thetag{4.4},
and let $A (\delta_\pm)$, $\Lambda (\delta_\pm)$, and $\widehat
N$ be the
operators acting on the Fock vectors by the rule
$$
\align
(\Phi,\Lambda(\delta_\pm)\Psi)
&=\lim_{\varepsilon\to + 0}
\sum_1^\infty\frac1{n!}\sum_{k=1}^n
\int_{(\R\setminus\{0\})^{n-1}}
\prod_{m\ne k} d\tau_m (\widetilde\Phi_n,\widetilde\Psi_n)_{\cal H}
\big|_{\tau_k=\pm\varepsilon},
\\
\FF_{\omega\to\tau} \bigl (A (\delta_\pm) \Psi\bigr)_n (\tau)
&=\lim_{\varepsilon\to+0}
\sum_{k=1}^n\widetilde{\Psi}_{n+1}(\tau)
\big|_{\tau_k=\pm\varepsilon}, \qquad \widehat N\Psi_n (\omega)
=n\Psi_n (\omega).
\endalign
$$
In this notation, the boundary condition~\thetag{4.4} in
$\Gamma^S$ acquires the form
$$
(\widehat N + 1 ^{-1} \bigl (I\otimes A (\delta_+) -W\otimes A (
\delta_-) \bigr) \Psi =L\otimes I\,\Psi.
\tag4.5
$$
This condition is imposed on  the Fock vector components
$\Psi_1,\Psi_2,\dots$
and extends condition \thetag{2.2} to the Fock space.

Let us prove that the operator
$$
\widehat H =iG\otimes I + I\otimes\widehat E + iL^* W\otimes A (\delta_-),
\qquad \widehat E =\int_{\R\setminus\{0\}} d\tau\, a^+ (\tau) a (
\tau) i\partial_\tau,
\tag4.6
$$
is symmetric in $\cal D_{W, L} $.
Let $\Phi,\Psi\in\cal D_{W, L} $, and let $B$ be a Hermitian
operator such that $\dom B\otimes I\supseteq\cal D_{W, L} $.
Integration by parts yields the following identity, where the 
integrated terms are 
expressed via the operators $\Lambda (\delta_\pm) $:
$$
\align
(\Phi, B\otimes\widehat E\Psi) - (B\otimes\widehat E\Phi,\Psi) =i\bigl (
\Phi, B\otimes\bigl (\Lambda (\delta_-) -
\Lambda(\delta_+)\bigr)\Psi\bigr)
\\\qquad
=i\sum_1^\infty\frac1{n!}\sum_{k=1}^n
\int_{(\R\setminus\{0\})^{n-1}}
\prod_{m\ne k} d\tau_m \bigl (\widetilde\Phi_n (\tau), B\widetilde\Psi_n
(\tau) \bigr)_{\cal H}
\big|_{\tau_k=+0}^{\tau_k=-0}.
\tag4.7
\endalign
$$
Using the boundary condition ~\thetag{4.5} for
$\widetilde\phi_n$ and $\widetilde\psi_n$, we find the
integrated terms in~\thetag{4.7}:
$$
\align
i\bigl (\Phi, B\otimes\bigl (\Lambda (\delta_+)
-\Lambda(\delta_-)\bigr)\Psi\bigr)
=i(\Phi,(W^*BW-B)\Lambda(\delta_)-\Psi)
\\\qquad+i(L\Phi,BL\Psi)
+ I (WA (\delta_-) \Phi, BL\Psi) + i (L\Phi, BWA (\delta_-) \Psi).
\tag4.8
\endalign
$$
In particular, Eq.~\thetag{4.8} is simplified for $B=I$ as
follows: 
$$
I\bigl (\Phi, I\otimes\bigl (\Lambda (\delta_+)
-\Lambda(\delta_-)\bigr)\Psi\bigr)
=i(\Phi,L^*L\Psi)
- (iL^* W\otimes A (\delta_-) \Phi,\Psi) + (\Phi, iL^* W\otimes A (
\delta_-) \Psi).
$$

Since $iG-iL^* L= (iG) ^* $, we obtain the following
identity, which means that the
operator $\widehat H$ is symmetric in $\cal D_{W, L} $:
$$
\align
(\Phi,\widehat H\Psi)
&=\bigl((I\otimes\widehat E \Phi,\Psi\bigr) + \bigl (\Phi,\{iG\otimes I
+ iL^* W\otimes A (\delta_-) \} \Psi\bigr) \\
&\qquad \qquad \qquad \qquad -i\bigl (\Phi,I\otimes\bigl (\Lambda (\delta_+)
-\Lambda(\delta_-)\bigr)\Psi\bigr)
\\
&=\bigl((I\otimes\widehat E \Phi,\Psi\bigr) + \bigl (\Phi, (iG) ^*
\otimes I\Psi\bigr) + (iL^* W\otimes A (\delta_-) \Phi,\Psi) = (
\widehat H\Phi,\Psi).
\endalign
$$

Let us find how the generator of the group $U_t$ acts on the range
of the resolvent. Set $\Psi\in\bold h$ and $\Phi=R_\mu h\otimes\psi (v) $.
By the definition of the generator, we have
$$
\alignat1
(\Psi,\widehat H\Phi)
&=\lim_{s\to+0}
\frac1 i\frac d {ds} (\Psi, U_s\Phi)
=\frac1 i\int_0^\infty dt\,
\sum_{n=0}^\infty\frac1{n!}
\int_{(\R\setminus\{0\})^n}d\tau\\
&\qquad\times \Bigl (\widetilde\psi_n (\tau),\frac d {ds}
e^{-(G+\mu)t-Gs-iL^*W\int_0^{t+s}
\widetilde v (-\tau) \, d\tau} \widetilde\phi_{n, t + s} (\tau) \Bigr
)_{\cal H} \Big|_{s=0}.
\tag4.9
\endalignat
$$
Note that the functions $\widetilde\phi_{n, t} (\tau) $ depend
on the differences
$\tau_k-t$.
Therefore,
$$
\frac d {dt} \widetilde\phi_{n, t} (\tau) =-\sum_{k=1}
^n\frac\partial {\partial\tau_k} \widetilde\phi_{n, t} (\tau)
=i\widehat E\widetilde\phi_{n, t} (\tau).
\tag4.10
$$
On the other hand,
$$
\widetilde{A(\delta_-)\phi_{n,t}}(\tau)
=n\widetilde v (-t) \phi_{n-1, t} (\tau)
\tag4.11
$$
by the definition of $A (\delta_-) $.
Now, taking into account definition~\thetag{4.9} and identities~\thetag
{4.10} and~\thetag{4.11}, we obtain
$$
\align
(\Psi,\widehat H\Phi)_{\bold h} &= \int_0^\infty dt \biggl (
\psi_0, iGe^{- (G + \mu) t-iL^* W\int_0^t \widetilde v (-\tau) \,
d\tau} h\biggr)_{\cal H}
+ \int_0^\infty dt
\sum_{n=1}^\infty\frac1{n!}
\int_{(\R\setminus\{0\}) ^n} d\tau \\
&\quad\qquad\times\biggl(\widetilde\psi_n(\tau),
e^{-(G+\mu)t-iL^*W\int_0^t
\widetilde v (-\tau) \, d\tau} \biggl (iG + iL^* W\widetilde v (-t)
+ I\sum_{k=1} ^n\frac\partial {\partial\tau_k} \biggr) \widetilde\phi_{
n, t} (\tau) \biggr)_{\cal H} \\
&=\bigl (\Psi,\{iG + iL^*
W\otimes A (\delta_-) + I\otimes\widehat E\} \Phi\bigr)_{\bold h};
\endalign
$$
that is, the generator $\widehat H$ of the group $U_t$
has the form \thetag{4.6}. Thus, we have proved the following theorem.

\proclaim{Theorem 4.1}  The symmetric operator
$$
\widehat H = iG\otimes I + I\otimes\widehat E + iL^* W\otimes A (\delta_-)
$$
in $\cal D_{W, L}$, where $G=iH + \frac12 L^* L$,
$H=\frac14 L^* KL-H_0$, is the generator of the one-parameter
unitary group~$U_t$.
\endproclaim

We must point out that the verification of the symmetry property
did not rely on the assumption  that $L$, $G$, and $ W$
commute and  can readily be extended to operators of the form
$$
\widehat H=iG + I\otimes\widehat E + i\sum_{l, m} L^*_l W_{l, m}
\otimes A_m (\delta_-)
\tag4.12
$$
with the boundary condition
$$
(\widehat N + 1) ^{-1} \biggl (I\otimes A_l (\delta_+) -\sum_mW_{l,
m} \otimes A_m (\delta_-) \biggr) \Psi =L_l\otimes I\Psi,
\tag4.13
$$
where $W=\{W_{l, m} \} $ is an
$ M\times M $-matrix  with entries in $\cal B (\cal H) $ such
that $W^* W=I$, and $\{A_m (g): g\in L_2 (\R),1\le m\le M\} $ are
the annihilation operators in $\Gamma^S(L_2(\R^M))$, which
commute for different $l$.

\head\bS5. The Markov evolution equation \endhead

In conclusion, we describe a derivation of the Markov
evolution equation from
the Cauchy  problem for  the Schr\"odinger equation
$$
\frac d {dt} \Psi (t) =\biggl (-G + iI\otimes\widehat E
-\sum_{l,m}L^*_lW_{l,m}
\otimes A_m (\delta_-) \biggr) \Psi (t)
$$
with the boundary conditions \thetag{4.13}.
Let $B\in\cal B (\cal H) $ be a Hermitian  operator,  and let $h, g\in D$.
Consider the equation 
$$
(G, P_t (B) h)_{\cal H} =\bigl (U_tg\otimes\Psi (0)| B\otimes I|
U_th\otimes\Psi (0) \bigr)_{\bold h}
$$
for the mean value.
From~\thetag{4.12} we have
$$
\align
\frac d {dt} (g, P_t (B) h)_{\cal H}
&=-\biggl(\biggl(G
+ \sum_{l, m} L^*_lW_{l, m} A_m (\delta_-) \biggr) U_t g\otimes\Psi
(0)| B\otimes I| U_th\otimes\Psi (0) \biggr) \\
&\qquad -\biggl (U_tg\otimes\Psi (0)| B\otimes I|
\biggl (G + \sum_{l, m} L^*_lW_{l,
m} A_m (\delta_-) \biggr)
U_th\otimes\Psi(0)\biggr)
\\
&\qquad+\biggl(U_tg\otimes\Psi(0)\bigl|
B\otimes\bigl (\Lambda (\delta_+) -\Lambda (\delta_-) \bigr) \bigr|
U_th\otimes\Psi(0))\biggr).
\tag5.1
\endalign
$$
Now we can use identity~\thetag{4.8}, which  in this
case can be  rewritten as follows:
$$
\multline
\biggl(\Phi, B\otimes\bigl (\Lambda (\delta_+) -
\Lambda(\delta_-)\bigr)\Psi\biggr)
=\sum_{l, m} \biggl (\Phi,\bigl (W^*_{l, m} BW_{l, m} -B\bigr)
\Lambda_m(\delta_-)\Psi\biggr)
\\
+\sum_l(L_l\Phi,BL_l\Psi)
+\sum_{l,m}\biggl(
(W_{l,m}A_m(\delta_-)\Phi,BL_l\Psi)
+ (L_l\Phi, BW_{l, m} A_m (\delta_-) \Psi) \biggr).
\endmultline
\tag5.2
$$
Taking into account the fact that $A (\delta_-) \Psi (0) =0$ and
$\Lambda(\delta_-)\Psi(0)=0$,
>from Eqs.~\thetag{5.1} and~\thetag{5.2} we obtain
$$
\frac d {dt} (g, P_t (B) h)_{\cal H} \Big|_{t=0} = (g,\cal L (B) h)
=-(Gg,Bh)-(g,BGh)+\sum_l(L_lg,BL_lg).
$$
Thus, we have obtained the generator $\cal L (\,\cdot\,) $ of
the Markov evolution equation  in Lindblad  form:
$$
\frac d {dt} P_t (B) =\cal L\bigl (P_t (B) \bigr), \quad \cal L (B)
=-G^* B-BG + \frac12\sum_lL^*_lBL_l, \quad
G=iH+\frac12\sum_lL^*_lL_l.
$$

This research was partially supported
by the Russian Foundation for Basic Research
under grant~No.~95-01-00784
  
\Refs

\item{1.}
G.~Lindblad,
``On the generators of quantum dynamical semigroups,''
{\sl Comm. Math. Phys.},
{\bf 48},
No.~2,
119--130
(1976).

\item{2.}
R.~L.~Hudson and K.~R.~Parthasarathy,
``Quantum Ito's formula and stochastic evolutions,''
{\sl Comm. Math. Phys.},
{\bf 93},
No.~3,
301--323
(1984).
  
\item{3.}
K.~R.~Parthasarathy,
{\sl An Introduction to Quantum Stochastic Calculus},
Birkhauser,
Basel
(1992).
  
\item{4.}
P.~A.~Meyer,
{\sl Quantum Probability for Probabilists},
Vol.~1338,
Lecture Notes in Math.,
Springer-Verlag,
Berlin
(1993).
  
\item{5.}
E.~B.~Davies,
{\sl Quantum Theory of Open Systems},
Acad. Press,
London
(1976).
  
\item{6.}
V.~Gorini, A.~Kossakovsky, and E.~C.~G.~Sudarshan,
``Completely positive dynamical semigroups of $n$-level systems,''
{\sl J. Math. Phys.},
{\bf 17},
No.~3,
821--825
(1976).
  
\item{7.}
C.~W.~Gardiner and M.~J.~Collett,
``Input and output in damped quantum systems: quantum statistical
differential equations and the master equation,''
{\sl Phys. Rev.~A},
{\bf 31},
3761--3774
(1985).
  
\item{8.}
{\sl Quantum Stochastic Processes and Open Systems}
(A.~N.~Kolmogorov and S.~P.~Novikov, editors)
[in Russian],
Vol.~42,
Ser. {\sl Mathematics},
Mir,
Moscow
(1988).

\item{9.}
R.~Z.~Khas{\mz}minskii, 
``Ergodic properties of recurrent diffusions and stabilization
of the Cauchy problem for a parabolic equation,''
{\sl Teor. Veroyatnost. i Primenen.} [{\sl Theory Probab. Appl.}],
{\bf 5},
No.~1,
196--214
(1960).
  
\item{10.}
K.~Ichihara,
{\sl Explosion Problems for Symmetric Diffusion Processes},
Vol.~1203,
Lecture Notes in Math.,
Springer-Verlag,
Berlin 
(1986),
pp.~75--89.

\item{11.} 
A.~M.~Chebotarev,
``Necessary and  sufficient conditions for conservativeness
of dynamical semigroup,''
in: Journal of Soviet Mathematics,
Vol. ~56,
No. ~5
(1991),
pp.~2697 - 2719.

\item{12.}
A.~M.~Chebotarev,
``Sufficient conditions of the  conservatism
of a minimal dynamical semigroup,''
{\sl Mat. Zametki\/} [{\sl Math. Notes\/}],
Vol. 52,
No.~4,
112--122.
(1992)

\item{13.}
A.~M.~Chebotarev, F.~Fagnola and A.~Frigerio,
``Towards a stochastic Stone's theorem,''
in: {\sl Stochastic Partial Differential Equations and Applications},
Vol.~268,
Pitman Res. Notes Math. Ser,
Longman Sci. Tech.,
Harlow
(1992)
pp.~86--97.

\item{14.}
A.~M.~Chebotarev and F.~Fagnola,
``Sufficient  conditions for conservativity of quantum dynamical
semigroups,''
{\sl J. Funct. Anal.},
{\bf 113},
No.~1,
131--153
(1993).

\item{15.}
A.~S.~Holevo,
``On conservativity of covariant dynamical semigroups,''
{\sl Rep. Math. Phys.},
{\bf 33},
95--100
(1993).

\item{16.}
B.~V.~Bhat \,R. and K.~R.~Parthasarathy,
``Markov dilations of non-conservative dynamical semigroups and a quantum
boundary theory,''
{\sl Ann. Inst. H.~Poincar\'e. Probab. Statist.},
{\bf 31},
No.~4,
601--651
(1995).

\item{17.}
A.~S.~Holevo,
``On the structure of covariant dynamical semigroups,''
{\sl J. Funct. Anal.},
{\bf 131},
255--278
(1995).

\item{18.}
A.~M.~Chebotarev, J.~K.~Garcia and R.~B.~Quezada,
``On the Lindblad equation with unbounded time-dependent coefficients'',
{\sl Mat. Zametki\/} [{\sl Math. Notes\/}],
{\bf 61},
No.~1,
125--140
(1997).

\item{19.}
A.~M.~Chebotarev,
``Minimal solutions in classical and quantum probability,''
in: {\sl Quantum Probability and Related Topics}, VII
(L.~ Accardi, editor),
World Scientific,
Singapore
(1992),
pp.~79--91.

\item{20.}
F.~Fagnola,
``Characterization of isometric and unitary weakly differentiable
cocycles in Fock space,''
in: {\sl Quantum Probability and Related Topics}, VIII,
Preprint No.~358,
UTM,
Trento
(1993),
pp.~143--164.

\item{21.}
B.~V.~Bhat \,R., F.~Fagnola and K.~B.~Sinha,
``On quantum extensions of semigroups of Brownian motions on a
half-line,''
{\sl Russian J. Math. Phys.},
{\bf 4},
No.~1,
13--28
(1996).

\item{22.}
J.~L.~Journ\'e,
``Structure des cocycles markoviens sur l'espace de Fock,''
{\sl Probab. Theory Related Fields},
{\bf 75},
291--316
(1987).

\item{23.}
A.~M.~Chebotarev,
``Symmetric form of the Hudson--Parthasarathy stochastic equation,''
{\sl Mat. Zametki\/} [{\sl Math. Notes\/}],
{\bf 60},
No.~5,
726--750
(1996).

\item{24.}
V.~D.~Koshmanenko,
``Perturbation of self-adjoint operators by singular bilinear forms,''
{\sl Ukrain. Mat. Zh.\/} [{\sl Ukrainian Math. J.\/}],
{\bf 41},
No.~1,
3--18
(1989).

\item{25.}
V.~D.~Koshmanenko,
{\sl Singular Bilinear Forms in
Perturbation Theory of Self-Adjoint Operators\/}
[in Russian],
Naukova Dumka,
Kiev
(1993).

\item{26.}
S.~Albeverio, W.~Karwowski and V.~Koshmanenko,
``Square powers of singularly perturbed operators,''
{\sl Math. Nachr.},
{\bf 173},
5--24
(1995).

\item{27.}
T.~Kato,
Perturbation   Theory   for   Linear   Operators,
Springer,
Heidelberg
(1976)

\item{28.}
F.~A.~Berezin,
{\sl The Method of Second Quantization\/}
[in Russian],
Nauka,
Moscow
(1986).

\endRefs

\authoraddress{M.~V.~Lomonosov Moscow State University, Quantum Statistics
Department, Moscow 119899}

\transl{author} 
\enddocument